\documentclass{scrartcl}

\usepackage{xcolor}
\colorlet{RefColor}{green!50!black}
\colorlet{LinkColor}{red!50!black}
\usepackage{hyperref}
\hypersetup{
  colorlinks = true,
  linkbordercolor = {white},
  linkcolor = LinkColor,
  anchorcolor=LinkColor,
  citecolor= RefColor,
  filecolor=cyan,
  menucolor=red,
  runcolor=cyan,
  urlcolor = LinkColor,
}
\usepackage[top=1in]{geometry}

\usepackage[numbers]{natbib}
\usepackage{authblk}
\usepackage{graphicx}%
\usepackage{multirow}%
\usepackage{amsmath,amssymb,amsfonts}%
\usepackage{amsthm}%
\usepackage{algorithm}%
\usepackage{algorithmicx}%
\usepackage{algpseudocode}%

\definecolor{mdred}{HTML}{a90033}
\definecolor{mdg}{HTML}{00aa23}

\newcommand{\R}{\mathbb{R}}

\newcommand{\bfs}{\boldsymbol{s}}

\newcommand{\reduced}[1]{\widetilde{#1}}

\newcommand{\nfull}{n}
\newcommand{\nsnapshots}{k}
\newcommand{\nred}{r}

\newcommand{\nredmod}{p}

\newcommand{\nconsider}{q}

\newcommand{\decoder}{g}
\newcommand{\encoder}{f}

\newcommand{\featuremap}{h}

\newcommand{\vecsymb}[1]{\boldsymbol{#1}}

\newcommand{\fullstate}{\vecsymb{s}}

\newcommand{\fullstatei}[1]{\vecsymb{s}^{(#1)}}

\newcommand{\redstate}{\widetilde{\vecsymb{s}}}
\newcommand{\redstatei}[1]{\widetilde{\vecsymb{s}}^{(#1)}}

\newcommand{\projection}[1]{\vecsymb{P}_{#1}}
\newcommand{\snapshots}{\vecsymb{S}}

\newcommand{\leftsings}{\vecsymb{\Phi}}
\newcommand{\singvals}{\vecsymb{\Sigma}}
\newcommand{\rightsings}{\vecsymb{\Psi}}

\newcommand{\leftsing}[1]{\vecsymb{\phi}^{(#1)}}
\newcommand{\scaledleftsing}[1]{\vecsymb{\phi}_{\sigma^{-1}}^{(#1)}}

\newcommand{\singval}[1]{\sigma_{#1}}

\newcommand{\Vlin}{\vecsymb{V}}
\newcommand{\scaledVlin}{\boldsymbol{V}^{\sigma^{-1}}}
\newcommand{\rescaledVlin}{\boldsymbol{V}^{\sigma}}

\newcommand{\Vnonlin}{\vecsymb{W}}

\newcommand{\indexin}[1]{\mathcal{I}_{#1}}
\newcommand{\indexout}[1]{\breve{\mathcal{I}}_{#1}}

\newcommand{\testvec}{\vecsymb{v}}

\newcommand{\norm}[1]{{\left\| #1 \right\|}}
\newcommand{\frobsq}[1]{\norm{#1}_\mathsf{F}^2}
\newcommand{\frob}[1]{\norm{#1}_\mathsf{F}}

\newcommand{\bfA}{\boldsymbol{A}}

\newcommand{\bfH}{\boldsymbol{H}}

\newcommand{\recloss}{\mathcal{J}_{\mathrm{rec}}}
\newcommand{\newloss}{\mathcal{J}_{\mathrm{tot}}}

\title{Operator Inference Aware Quadratic Manifolds with Isotropic Reduced Coordinates for Nonintrusive Model Reduction}

\author[1]{Paul Schwerdtner\thanks{paul.schwerdtner@nyu.edu}}
\author[2]{Prakash Mohan\thanks{prakash.mohan@nrel.gov}}
\author[2]{Julie Bessac\thanks{julie.bessac@nrel.gov}}
\author[2]{Marc T. Henry de Frahan\thanks{marc.henrydefrahan@nrel.gov}}
\author[1]{Benjamin Peherstorfer\thanks{pehersto@cims.nyu.edu}}

\affil[1]{Courant Institute of Mathematical Sciences, New York University, 251 Mercer Street, New York, NY 10012, USA}
\affil[2]{Computational Science Center, National Renewable Energy Laboratory, 15013 Denver West Parkway, Golden, CO 80401, USA}

\date{July 2025}

\begin{document}

\maketitle
\begin{abstract}
  \noindent
Quadratic manifolds for nonintrusive reduced modeling are typically trained to minimize the reconstruction error on snapshot data, which means that the error of models fitted to the embedded data in downstream learning steps is ignored. In contrast, we propose a greedy training procedure that takes into account both the reconstruction error on the snapshot data and the prediction error of reduced models fitted to the data. Because our procedure learns quadratic manifolds with the objective of achieving accurate reduced models, it avoids oscillatory and other non-smooth embeddings that can hinder learning accurate reduced models. 
Numerical experiments on transport and turbulent flow problems show that quadratic manifolds trained with the proposed greedy approach lead to  reduced models with up to two orders of magnitude higher accuracy than quadratic manifolds trained with respect to the reconstruction error alone.  
\end{abstract}
\vspace{0.5cm}

\vspace{0.5cm}
\newpage

\section{Introduction}
Learning reduced models from data in a nonintrusive fashion is an important problem in science and engineering \cite{Ghattas_Willcox_2021,doi:10.1137/1.9781611976083,annurev:/content/journals/10.1146/annurev-fluid-121021-025220}. A typical approach is to first learn an encoder-decoder pair, embed the training snapshot trajectories  with the learned encoder, and then fit a reduced dynamical-system model to the embedded trajectories.
However, the decomposition of the training process into first learning an encoder-decoder pair for the embedding and only sub-sequentially learning a model of the dynamics typically means that the encoder-decoder pair are trained with the objective of accurately approximating the training data, rather than taking the reduced-model prediction error into account.
Thus, the encoder-decoder pair can overfit to achieving a low reconstruction error on the training data by learning embeddings of the snapshot trajectories that are non-smooth, which means that learning a reduced model can become challenging. Correspondingly, it has been observed that learning embeddings and models together can be beneficial; see, e.g., \cite{NIPS2017_3a835d32,Lusch2018,doi:10.1137/18M1177846,doi:10.1073/pnas.1906995116}. In the context of intrusive model reduction with linear approximations, there is work that optimizes the reduced basis with respect to the model prediction error \cite{DBLP:conf/amcc/Borggaard06},  quantities of interest \cite{BUITHANH2007880}, and to achieve stability \cite{BARONE20091932}; however, we focus here on the setting of nonintrusive model reduction and nonlinear approximations.

In this work, we consider encoder-decoder pairs corresponding to quadratic manifolds and propose a greedy training approach that takes the reconstruction error as well as the model prediction error into account. We build on the greedy approach in \cite{SchwerdtnerP2024Greedy} and propose a new selection criterion that includes the prediction error of a fitted operator inference model \cite{PeherstorferW2016Data-driven}. The new selection criterion nudges the greedy approach to find an encoder-decoder pair that leads to embedded trajectories to which operator inference models can be fitted well. Numerical experiments with transport and turbulent flow problems demonstrate that the proposed operator inference-aware greedy approach can lead to almost two orders of magnitude higher accuracy in operator inference models compared to fitting quadratic manifolds based on the reconstruction error alone.

Nonintrusive model reduction has emerged from intrusive, projection-based model reduction~\cite{Antoulas2005Approximation,RozzaHP2008Reduced,BennerGW2015survey,HesthavenRS2016Certified} and typically requires only snapshot data to construct reduced models. Examples for nonintrusive model reduction methods include dynamic mode decomposition (DMD)~\cite{rowley2009spectral,Schmid2010Dynamic,KutzBBP2016Dynamic,brunton2016koopman}, the Loewner framework~\cite{AntoulasLI2017tutorial,MayoA2007framework,AntoulasBG2020Interpolatory}, optimization-based model reduction \cite{SchwerdtnerV2023SOBMOR,optimizationparametric}, and data-driven balanced truncation approaches \cite{doi:10.1137/21M1411081,reiter2023generalizationsdatadrivenbalancingsample}; see \cite{Ghattas_Willcox_2021,annurev:/content/journals/10.1146/annurev-fluid-121021-025220} for more comprehensive overviews on nonintrusive model reduction. We focus on operator inference for nonintrusive model reduction, which was introduced in~\cite{PeherstorferW2016Data-driven} and has become a building block for a series of works; see, e.g.,  \cite{McQuarrie03042021,benner2020deim,QianKPW2020LiftLearn,GUO2022115336,doi:10.1137/21M1452810,Peherstorfer2020Sampling,benner2022incompressible,doi:10.1098/rsta.2021.0206,doi:10.1137/21M1439729,UyHP2023Operator}. In particular, progress has been made on imposing and preserving structure in operator inference models~\cite{SharmaWK2022Hamiltonian,GruberT2023Canonical,SHARMA2024134128,SAWANT2023115836}.
Operator inference has also served as a starting point for constructing nonintrusive reduced models for a range of applications from large-scale combustion-dynamics~\cite{McQuarrie03042021,doi:10.2514/1.J063715,FARCAS2025109619} to plasma physics~\cite{10.1063/5.0225584} to shallow water wave problems~\cite{yildiz2021shallow}.
A major step forward is the work \cite{GeelenWW2023Operator} that develops the first operator inference approach based on nonlinear approximations with quadratic manifolds, which allows reducing transport-dominated problems; see \cite{P22AMS} for a brief survey on intrusive model reduction with nonlinear approximations. Quadratic manifolds have been introduced to model reduction in \cite{JainTRR2017quadratic,RutzmoserRTJ2017Generalization}. A range of nonlinear model reduction methods that build on quadratic or similar forms of nonlinear decoders has appeared since then~\cite{BarnettF2022Quadratic,CRMECA_2023__351_S1_357_0,GeelenBWW2024Learning,GeelenBW2023Learning,BennerGHP-D2023quadratic,SchwerdtnerP2024Greedy,SharmaMBGGK2023Symplectic,SchwerdtnerGP2024Empirical,WederSP2024Nonlinear,ballout:hal-05143191,diaz2025interpretableflexiblenonintrusivereducedorder}. 
In all of these works, quadratic manifolds are constructed based on the reconstruction error alone. In contrast, we build on the greedy approach introduced in \cite{SchwerdtnerP2024Greedy} and introduce a new greedy criterion that also takes the prediction error of an operator inference model into account. Thus, the encoder-decoder pair is trained based on the reconstruction as well as the model prediction error.  

This manuscript is organized as follows. We first discuss preliminaries and provide the problem formulation in Section~\ref{sec:Prelim}. The novel greedy approach is introduced in Section~\ref{sec:OpInfGreedy} and consists of a step to obtain isotropic reduced coordinates and a new selection criterion that takes the operator-inference model prediction error into account. Numerical experiments on a transport and a turbulent flow problem are provided in Section~\ref{sec:NumExp}, where the proposed approach achieves orders of magnitude higher accuracy that the regular greedy approach from \cite{SchwerdtnerP2024Greedy} as well as other approaches for constructing quadratic manifolds. Conclusions are drawn in Section~\ref{sec:Conc}

\section{Preliminaries}\label{sec:Prelim}

\subsection{Data reduction with quadratic decoder functions}
Let $\fullstate \in \R^{\nfull}$ be a high-dimensional state vector of the full model, i.e., of the model of which we would like to construct a reduced model.
Furthermore, let  $\encoder: \R^{\nfull} \to \R^{\nred}$ be an encoder that maps a state $\fullstate$ onto $\R^{\nred}$ with $\nred \ll \nfull$. A corresponding decoder is denoted as $\decoder: \R^{\nred}\to\R^{\nfull}$.
A standard approach is   training encoder-decoder pairs with respect to the reconstruction error 
\begin{align}
  \label{eq:rec_err}
  e(\fullstate)={\|\decoder(\encoder(\fullstate))-\fullstate\|}_2^2
\end{align}
for states that are typically attained during simulations of the full model.
For example, linear encoder and decoder functions can be obtained via the singular value decomposition (SVD) of the snapshot matrix $\snapshots = [\fullstatei{1}, \dots, \fullstatei{\nsnapshots}]$, which contains the states of the full model over a time trajectory. 
Let  $\leftsings \singvals \rightsings^\top=\snapshots$ be the reduced SVD of $\snapshots$ where $\leftsings \in \mathbb{R}^{\nfull \times r_{\text{max}}}$ is the matrix of left singular vectors, $\singvals \in \mathbb{R}^{r_{\text{max}} \times r_{\text{max}}}$ is the diagonal matrix of singular values $\sigma_1 \geq \cdots \geq \sigma_{r_{\text{max}}} > 0$, and $\rightsings \in \mathbb{R}^{k \times r_{\text{max}}}$ is the matrix of right singular vectors.
The linear encoder and decoder functions that minimize the averaged reconstruction error (over all linear encoder and decoder functions) over the training data are given by $\encoder_{\Vlin}(\fullstate)=\Vlin^\top \fullstate$ and $\decoder_{\Vlin}(\redstate)=\Vlin\redstate$.

As shown in the context of model reduction in~\cite{GeelenWW2023Operator,JainTRR2017quadratic,RutzmoserRTJ2017Generalization, GeelenBWW2024Learning, GeelenBW2023Learning, BarnettF2022Quadratic, SchwerdtnerP2024Greedy}, augmenting the decoder function with an additive nonlinear correction term can achieve lower reconstruction errors than using just a linear decoder even when the encoder is kept linear; in particular, a detailed discussion that a linear encoder is sufficient is provided in \cite{SchwerdtnerP2024Greedy}.
In the following, we specifically consider nonlinear decoders based on nonlinear feature maps, which have the form
\begin{align}
  \label{eq:nonlinear_decoder}
  \decoder_{\Vlin, \Vnonlin}(\redstate)=\Vlin \redstate + \Vnonlin \featuremap(\redstate),
\end{align}
where $\featuremap: \R^{\nred}\to\R^{\nredmod}$ is a nonlinear feature map onto a $\nredmod$-dimensional feature space $\R^{\nredmod}$, $\Vnonlin \in \R^{\nfull \times \nredmod}$ is a weighting matrix, and $\Vlin \in \R^{\nfull \times \nred}$ is a matrix with orthonormal columns. We refer to $\Vlin$ as the basis matrix and to $\Vnonlin$ as the coefficient matrix.

Throughout this article, we follow \cite{GeelenWW2023Operator,BarnettF2022Quadratic,GeelenBWW2024Learning,GeelenBW2023Learning,JainTRR2017quadratic} and only use as feature map the quadratic function
\begin{align}
  \featuremap: \R^{\nred} \to \R^{\nred(\nred+1)/2}, \boldsymbol{x} \mapsto
  [x_1x_1,\,x_1x_2,\,\dots,\,x_1x_{\nred},\,x_2x_2,\,\dots, \,x_{\nred}x_{\nred}]^\top,
\end{align}
which is often called the condensed Kronecker product. The nonlinear decoder maps into the set 
\begin{equation}\label{eq:Prelim:ManifoldDef}
       \mathcal{M}_{\nred}(\Vlin, \Vnonlin)=\{ \Vlin \redstate + \Vnonlin \featuremap(\redstate) \,|\, \redstate \in \R^{\nred} \} \subset \R^{\nfull}\,,
\end{equation}
which includes the subspace spanned by the columns $\Vlin$ and additionally elements that can be described with linear combinations of the columns of the weight matrix $\Vnonlin$. Following convention, we refer to the set $\mathcal{M}_{\nred}(\Vlin, \Vnonlin)$ as quadratic manifold because the feature map $h$ is quadratic in the following \cite{JainTRR2017quadratic,GeelenWW2023Operator, BarnettF2022Quadratic}.
 
For a given basis matrix $\Vlin$, the coefficient matrix $\Vnonlin$ can be computed by solving the linear least-squares problem
    \begin{align}
      \label{eq:qm_lstsq_problem}
      \min\limits_{\Vnonlin \in \R^{\nfull \times \nredmod}} \frobsq{\projection{\Vlin}(\snapshots) + \Vnonlin \featuremap(\Vlin^\top \snapshots)-\snapshots} + \gamma \frobsq{\Vnonlin},
    \end{align}
where $\projection{\Vlin}$ denotes the orthogonal projection onto a subspace of $\R^{\nfull}$ with basis $\Vlin$, $\gamma > 0$ is a regularization parameter, and $\frob{\cdot}$ denotes the Frobenius norm. In~\eqref{eq:qm_lstsq_problem}, and throughout this article, $\featuremap$ is applied column-wise to its matrix argument.

\subsection[Greedily training basis matrix V]{Greedily training the basis matrix $\Vlin$}
One approach for training $\Vlin$ from snapshot data $\snapshots$ is selecting the leading $\nred$ left-singular vectors of $\snapshots$ as columns of $\Vlin$ \cite{GeelenWW2023Operator,BarnettF2022Quadratic}.
However, it has been shown in \cite{SchwerdtnerP2024Greedy}  that in many cases a lower reconstruction error on the snapshot data can be achieved by greedily selecting the columns of $\Vlin$ as follows:  
Let us consider the objective function
\begin{align}
  \label{eq:column_objective}
  \recloss(\testvec,\Vlin)=\min\limits_{\Vnonlin \in \R^{\nfull \times \nred(\nred+1)/2}}\frobsq{\projection{[\Vlin,\testvec]}(\snapshots) + \Vnonlin\featuremap(\encoder_{[\Vlin,\testvec]}(\snapshots))-\snapshots} + \gamma \frobsq{\Vnonlin}\,.
\end{align}
At iteration $i = 1, \dots, \nred$, the left-singular vector $\leftsing{j_i}$ out of the first $q > \nred$ vectors is selected that minimizes~\eqref{eq:column_objective},
\begin{align}
  \min\limits_{j_i=1,\dots,\nconsider}\recloss(\leftsing{j_i}, \Vlin_{i-1}).
\end{align}
The greedy algorithm proposed in \cite{SchwerdtnerP2024Greedy} starts at iteration $i=1$ with the empty matrix $\Vlin_0$. After $\nred$ iterations, the greedy basis matrix $\Vlin=[\leftsing{j_1}, \dots, \leftsing{j_{\nred}}]\in \R^{\nfull \times \nred}$ is obtained, which contains as columns the left-singular vectors of $\snapshots$ with indices $j_1, \dots, j_{\nred}$. Then the weight matrix $\Vnonlin$ is computed using~\eqref{eq:qm_lstsq_problem}.
Once the $\Vlin$ and $\Vnonlin$ have been obtained, the linear encoder $f_{\Vlin}$ and quadratic decoder $g_{\Vlin,\Vnonlin}$ are fully defined. We refer to \cite{SchwerdtnerP2024Greedy} for a detailed discussion of the greedy algorithm and for a comparison to selecting the leading $\nred$ left-singular vectors. 

\subsection{Operator inference for nonintrusive model reduction}\label{sec:Prelim:OpInf}
Following \cite{GeelenWW2023Operator}, a reduced model can be constructed based on linear encoder and quadratic decoder functions with operator inference \cite{PeherstorferW2016Data-driven}. 
Given are the snapshot data $\snapshots = [\fullstatei{1}, \dots, \fullstatei{\nsnapshots}]$ from the full model. 
Let $[\redstatei{1}, \dots, \redstatei{\nsnapshots}]\in\R^{\nred \times \nsnapshots}$ denote the encoded (reduced) states, which are obtained as 
\[
\redstatei{i} = f_{\Vlin}(\fullstatei{i})\,,\qquad i = 1, \dots, \nsnapshots\,.
\]
Furthermore, define the matrices  $\tilde{\snapshots}_{-} = [\redstatei{1}, \dots, \redstatei{\nsnapshots-1}]$ and $\tilde{\snapshots}_+ = [\redstatei{2}, \dots, \redstatei{\nsnapshots}]$. 
Building on the formulation of operator inference used in \cite{GeelenWW2023Operator}, one fits the matrices $\bfA \in \mathbb{R}^{\nred \times \nred}$ and $\bfH \in \mathbb{R}^{\nred \times (\nred-1)\nred/2}$ as

\begin{align}
  \label{eq:opinf_lstsq}
  \min\limits_{\bfA \in \R^{\nred\times\nred}, \bfH \in \R^{\nred \times \nredmod}}
  \frobsq{\tilde{\snapshots}_{+}-\bfA \tilde{\snapshots}_{-}-\bfH \featuremap(\tilde{\snapshots}_{-})}+\gamma_{\bfA} \frobsq{\bfA}+\gamma_{\bfH} \frobsq{\bfH},
\end{align}
where $\gamma_{\bfA}$ and $\gamma_{\bfH}$ are regularization parameters.
Additional terms such as constant terms can be included as well. However, to keep our exposition concise, we only use linear and quadratic terms.
Once a $\bfA$ and $\bfH$ are obtained, they give rise to a reduced model of the form
\[
\hat{\bfs}^{({j + 1})} = \bfA\hat{\bfs}^{(j)} + \bfH h(\hat{\bfs}^{(j)})\,,\qquad j = 0, \dots, k - 1\,,
\]
that allows predicting a trajectory $\hat{\snapshots} = [\hat{\bfs}^{(1)}, \dots, \hat{\bfs}^{(k)}]$ for a new initial condition $\hat{\bfs}^{(0)} = f_{\Vlin}(\fullstate^{(0)}_\textrm{unseen})$ corresponding to some unseen full initial state $\fullstate^{(0)}_\textrm{unseen} \in \R^{\nfull}$.

\subsection{Problem statement}\label{sec:Prelim:Problem}

\begin{figure}
  \begin{center}
    \includegraphics[width=0.99\columnwidth]{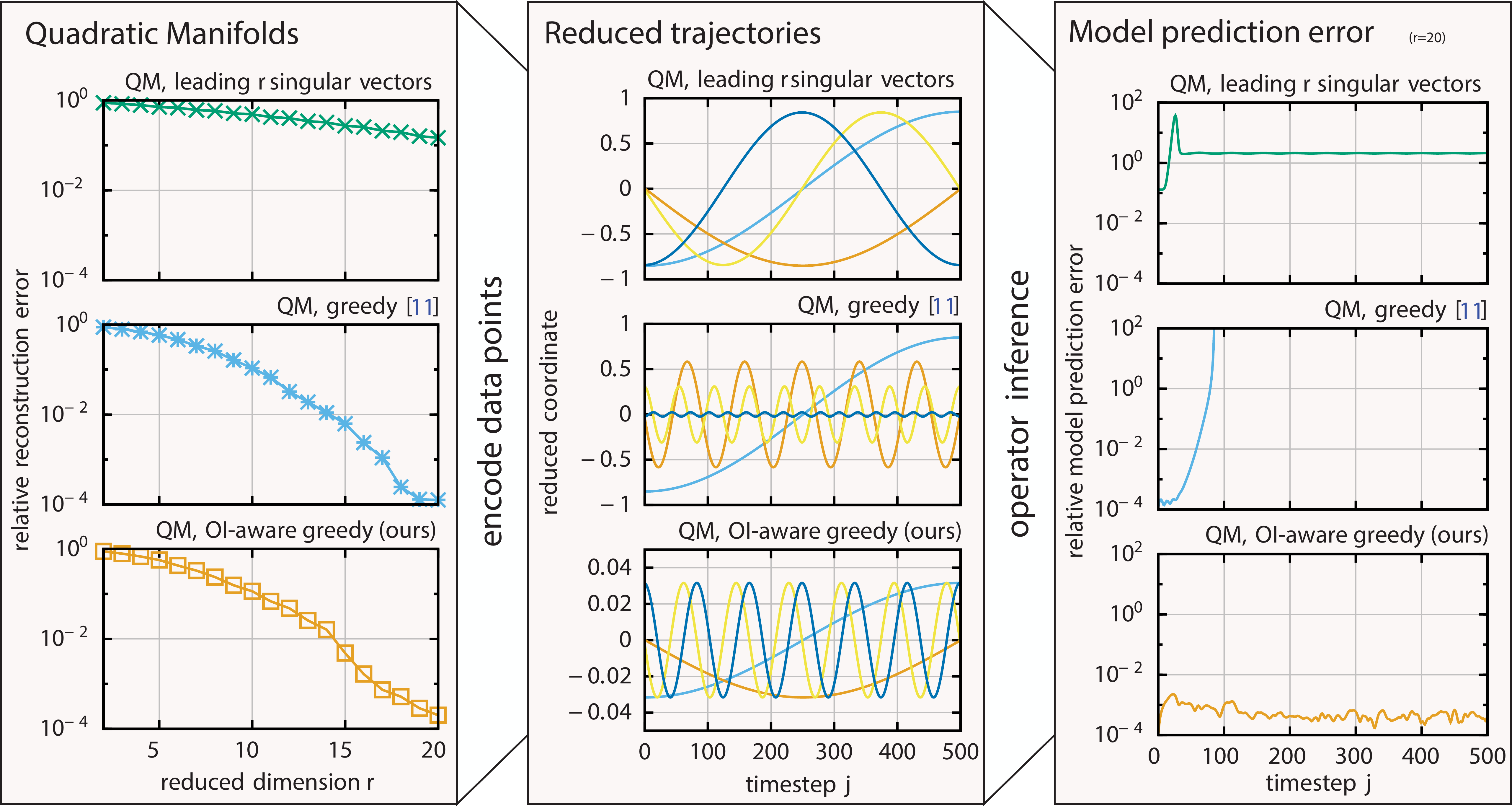}
  \end{center}
  \caption{Transport problem: The proposed operator inference-aware greedy approach (last row) considers both the reconstruction and prediction errors when training encoder-decoder pairs, which results in quadratic-manifold embeddings that are well-suited for nonintrusive model reduction with operator inference.}
  \label{fig:problem_statement}
\end{figure}

Using the greedy algorithm from \cite{SchwerdtnerP2024Greedy} to find $\Vlin$ can lead to an encoder-decoder pair that achieves a lower reconstruction error \eqref{eq:rec_err} than using the $\nred$ leading left-singular vectors only. However, even though the embedding is more accurate in terms of the reconstruction error, it can be more challenging to fit an operator inference model to the embedded data:
Figure~\ref{fig:problem_statement} shows that the coordinates of the reduced states obtained when choosing the first singular vectors as basis $\Vlin$ are similar in rate of change and magnitude such that reduced operators can be fitted well with operator inference. In contrast, when training the encoder-decoder pair with the greedy method that also selects singular vectors with higher indices into the matrix $\Vlin$, then the reduced coordinates corresponding to the embedded snapshots can oscillate more rapidly, which means it becomes more challenging to fit a reduced operator that can predict them.
In particular, because a first-order dynamical system is fitted to trajectories in operator inference, it cannot capture reduced coordinates that oscillate with starkly different frequencies in each dimension because at least two modes with similar frequencies are required so that a first-order system can capture the sine and cosine components required for oscillatory dynamics. Thus, focusing the training of the encoder-decoder pair of a quadratic manifold on the reconstruction error alone, can lead to difficulties when fitting an operator inference model downstream.

\section{Operator inference-aware greedy training of quadratic manifolds}\label{sec:OpInfGreedy}

We now introduce a new greedy selection criterion that balances reconstruction error and model prediction error.
We achieve this by fitting operator-inference models during the greedy iterations to compute the model prediction error.  
First, we introduce a modification to the encoder and decoder function that lead to isotropic reduced coordinates, which are  normalized in a sense that is suitable for operator inference with quadratic manifolds. Second, we define the new selection criterion and describe the corresponding greedy algorithm. 

\subsection{Isotropic reduced coordinates}\label{sec:OIGreedy:Isotropic}

When applying the linear encoder $\encoder_{\Vlin}$ to the snapshot data $\snapshots$, the reduced coordinates $\reduced{\snapshots}$ are anisotropic because the rows of $\reduced{\snapshots}$ are scaled by their corresponding singular value:  $\encoder_{\Vlin}(\snapshots)=\Vlin^\top \leftsings \singvals \rightsings^\top= \singvals_{[:\nred, :\nred]}\rightsings^\top_{[:\nred, :]}$. Notice that $\singvals_{[:\nred, :\nred]}$ denotes the upper left  $\nred \times \nred$ block of the matrix $\singvals$ and $\rightsings^\top_{[:\nred, :]}$ the first $\nred$ rows of $\rightsings^\top$. 
When applying operator inference using a linear decoder $\decoder_{\Vlin}$, the anisotropy due to the multiplication with the singular values can be useful because the scaling implicitly induces weights into the least-squares problem~\eqref{eq:opinf_lstsq} so that error corresponding to leading modes are overweight. 

When using quadratic decoders; however, the setting is different: The matrix  $\Vnonlin$ is the solution of a least-squares problem~\eqref{eq:qm_lstsq_problem} that does not in general consist of orthonormal columns. Moreover, the feature map $\featuremap$ combines reduced coordinates corresponding to leading and later left-singular vectors selected by the greedy method.
That is why for operator inference with quadratic decoder functions, the magnitude of the reduced coordinates is not necessarily indicative of the importance it carries for the reconstruction accuracy.
It is noteworthy that such a multi-scale effect did not arise in the earlier works such as~\cite{GeelenWW2023Operator}, where the leading singular vectors where used to form the basis matrix $\Vlin$, in contrast to the greedy method that also selects later left-singular vectors into $\Vlin$. In particular, if there is a slow singular value decay, as in the examples considered in \cite{GeelenWW2023Operator}, then the reduced coordinates corresponding to the first few leading singular vectors will be close to isotropic, which can also be seen in the upper row in Figure~\ref{fig:problem_statement}.
In contrast, the greedy column selection includes leading and later singular vectors and later later singular vectors can yield small-magnitude reduced coordinates even if the singular value decay slowly. But the corresponding reduced coordinates are still important for accuracy.

Because reduced coordinates with small magnitude can still be important for decoding with quadratic decoder, we introduce isotropic reduced coordinates as follows: We scale $\Vlin$ such that each row of $\reduced{\snapshots}$ has unit $\|\cdot\|_2$ norm. Notice that this is different from scaling and normalizing the snapshots in $\snapshots$, which is common in model reduction.
Isotropic reduced coordinates can be achieved by scaling each column of $\Vlin=[\leftsing{j_1}, \dots, \leftsing{j_{\nred}}]$ by the inverse of the corresponding singular value $\singval{j_1}, \dots, \singval{j_{\nred}}$.
We denote the scaled singular vectors by $\scaledleftsing{j}=\leftsing{j}/\singval{j}$ and the scaled basis by $\scaledVlin=[\scaledleftsing{j_1}, \dots, \scaledleftsing{j_{\nred}}]$.
Consequently, we denote the linear encoder that leads to isotropic reduced coordinates by $\encoder_{\scaledVlin}$ such that $\encoder_{\scaledVlin}\fullstate=(\scaledVlin)^\top\fullstate$.
For the decoders, we introduce $\rescaledVlin=[\singval{j_1}\leftsing{j_1}, \dots, \singval{j_{\nred}}\leftsing{j_{\nred}}]$. The linear and quadratic decoders that lift isotropic reduced coordinates are then given by
\[
  \decoder_{\rescaledVlin}(\redstate)=\rescaledVlin\redstate
\]
and
\[
  \decoder_{\rescaledVlin, \Vnonlin}(\redstate)=\rescaledVlin\redstate+ \Vnonlin \featuremap(\redstate),
\]
respectively.

\subsection{Balancing reconstruction error and model prediction error}
We now introduce a new greedy column selection objective that takes into account the model prediction error of an operator-inference model learned from the corresponding reduced trajectory; additionally to the reconstruction error accounted for by the objective \eqref{eq:column_objective}. 
Moreover, the new criterion is based on isotropic reduced coordinates that were introduced in the previous section.
Let $i = 1, \dots, \nred$ be the greedy iterations. 
At iteration $i$, we have the current $\scaledVlin_{i-1}$ and $\Vnonlin_{i-1}$ from the previous iteration, which give rise to the encoder $f_{\scaledVlin_{i-1}}$ and the corresponding decoder $g_{\rescaledVlin_{i-1}, \Vnonlin_{i-1}}$. 
We then use the current encoder $f_{\scaledVlin_{i-1}}$ to construct the reduced trajectory
\[
\reduced{\snapshots}_{i-1} = f_{\scaledVlin_{i-1}}(\snapshots) \in \R^{\nred \times \nsnapshots}
\]
from the given snapshots $\snapshots$. We then train an operator-inference model using the reduced trajectory $\reduced{\snapshots}_{i-1}$ following Section~\ref{sec:Prelim:OpInf} to obtain $\bfA_{i-1}$ and $\bfH_{i-1}$.
Then, the prediction of a model learned with operator inference is 
\begin{align}
  \label{eq:opinf_rollout}
  \hat{\redstate}^{(j+1)}=\bfA_{i-1} \hat{\redstate}^{(j)} + \bfH_{i-1} \featuremap(\hat{\redstate}^{(j)}), \qquad \hat{\redstate}^{(0)}=\redstate^{(0)}.
\end{align}
The corresponding trajectory is  $\hat{\reduced{\snapshots}} = [\hat{\redstate}^{(1)}, \dots, \hat{\redstate}^{(\nsnapshots)}]$.
We then use the following objective function to select the next left-singular vector to use as column $i$ in $\Vlin$, 
\begin{align}
  \label{eq:newloss}
  \newloss(\testvec, \scaledVlin) = \recloss(\tilde \testvec, \scaledVlin)+ \gamma_{\mathrm{op}} \frobsq{\reduced{\snapshots}-\hat{\reduced{\snapshots}}},
\end{align}
where $\reduced{\snapshots}=\encoder_{[\scaledVlin, \testvec]}(\snapshots)$ consists of isotropic reduced coordinates, when $\testvec$ is chosen from $[\scaledleftsing{1}, \dots, \scaledleftsing{\nsnapshots}]$. The regularization parameter $\gamma_{\mathrm{op}}$ can be selected based on strategies such as grid search and cross validation on a validation data set.
The objective \eqref{eq:newloss} takes the reconstruction error and the model prediction error into account when selecting the next column of $\scaledVlin$.
Note that~\eqref{eq:newloss} does not use the residual of the least-squares problem~\eqref{eq:opinf_lstsq} but the error of the trajectory $\hat{\reduced{\snapshots}}$.
This is because the residual is equivalent to the one-step prediction error, which is, in our setting, not informative of long term stability and accuracy of the learned operator inference model; see also \cite{UyHP2023Operator}. Finally, let us add that the approach can be extended to take into account the model prediction error over various different initial conditions. One way to achieve this extension is to evaluate the reduced model \eqref{eq:opinf_rollout} on several  initial conditions and then stack the corresponding trajectories to compute the error in the Frobenius norm in \eqref{eq:newloss} over all trajectories. 

\subsection{Speeding up the greedy objective function evaluation}
To determine a new column for the basis $\scaledVlin$, in iteration $i  = 1, \dots, \nred$, the objective function~\eqref{eq:newloss} is evaluated $\nconsider+1-i$ times.
The model prediction error can be computed at a cost that scales as $\mathcal{O}(\nred^3 \nsnapshots)$. Notice that the greedy approach to construct $\scaledVlin$ uses only evaluations of the objective function~\eqref{eq:newloss} and does not involve any gradient computations. This means there is no need to differentiate through the time stepping of the operator inference model \eqref{eq:opinf_rollout}. 

The costs for evaluating $\recloss$ scale as $\mathcal{O}(\nfull \nsnapshots^2)$, which can be substantial if the full state dimension $\nfull$ is high.
To reduce the costs, we use a similar insight as in~\cite{SchwerdtnerP2024Greedy}: Let $\indexin{}=\{j_1, \dots, j_{i}\}$ be the set of singular vectors included in $\scaledVlin$ and let $\indexout{}=\{1,\dots,\nsnapshots\} \setminus \indexin{}$ be its complement.
Note that the reconstruction error of the linear encoder-decoder pair can be written as
\begin{align}
  \decoder_{\rescaledVlin}(\encoder_{\scaledVlin}(\snapshots)) = \leftsings_{\indexout{}}\singvals_{\indexout{}}\rightsings_{\indexout{}}^\top.
\end{align}
Therefore, when $\recloss(\scaledleftsing{j'}, [\scaledleftsing{j_1}, \dots, \scaledleftsing{j_{i-1}}])$ with $\recloss$ defined in~\eqref{eq:column_objective} is evaluated, this leads to the least-squares problem
\begin{align}
  \label{eq:pre_quick_column}
  \min\limits_{\Vnonlin \in \R^{\nfull \times \nred(\nred+1)/2}} \frobsq{\leftsings_{\indexout{}}\singvals_{\indexout{}}\rightsings_{\indexout{}}^\top + \Vnonlin\featuremap(\rightsings_{\indexin{}}^\top)} + \gamma \frobsq{\Vnonlin},
\end{align}
where $\indexin{}=\{j_1, \dots, j_{i-1}, j'\}$. Here we used that $\encoder_{\scaledVlin}(\snapshots)=\rightsings_{\indexin{}}^\top$.
Considering that $\leftsings_{\indexout{}}$ has orthonormal columns and that it has more rows than columns, we realize that the minimizer of~\eqref{eq:pre_quick_column} is the same as the minimizer of
\begin{align}
  \label{eq:quick_column}
  \min\limits_{\Vnonlin' \in \R^{\nsnapshots-i \times \nred(\nred+1)/2}} \frobsq{\singvals_{\indexout{}}\rightsings_{\indexout{}}^\top + \Vnonlin'\featuremap(\rightsings_{\indexin{}}^\top)} + \gamma \frobsq{\Vnonlin'},
\end{align}
where $\Vnonlin'$ is the description of $\Vnonlin$ in the basis $\leftsings_{\indexout{}}$.
Thus, it suffices to evaluate~\eqref{eq:quick_column} when determining the next scaled column to add to the iteratively constructed bases $\scaledVlin$.
In this way, for estimating the reconstruction error, we can use the modified objective function
\begin{align}
    \recloss'(\scaledleftsing{j}, \scaledVlin_i) = 
    \min\limits_{\Vnonlin' \in \R^{\nfull \times \nred(\nred+1)/2}}
    \frobsq{
    \singvals_{\indexout{i}\setminus{j}}
    \rightsings_{\indexout{i}\setminus{j}}^\top
    +
    \Vnonlin' \featuremap(\rightsings^\top_{\indexin{i}\cup j})
    }
    + \gamma \frobsq{\Vnonlin'},
\end{align}
where $\scaledVlin_i=[\scaledleftsing{j_1}, \dots, \scaledleftsing{j_i}]$ with $\indexin{i}=\{j_1, \dots, j_i\}$ consists of the previously selected scaled singular vectors.
In combination with the operator-inference model prediction error, we obtain the objective
\begin{align}
\label{eq:newlossfast}
    \newloss'(\scaledleftsing{j}, \scaledVlin_i)=
    \recloss'(\scaledleftsing{j}, \scaledVlin_i) + 
    \gamma_{\mathrm{op}}
    \frobsq{\rightsings_{\indexin{i}}^\top-\hat{\reduced{\snapshots}}},
\end{align}
where $\hat{\reduced{\snapshots}}$ are the operator-inference predictions of the training trajectory $\encoder_{[ \scaledVlin_i,\scaledleftsing{j}]}(\snapshots)$.

\subsection{Algorithm}

Algorithm~\ref{alg:ourquadmani} summarizes the proposed approach. 
The algorithm takes as input the snapshot matrix $\snapshots$, the reduced dimension $\nred$, the regularization parameters $\gamma, \gamma_{\bfA}, \gamma_{\bfH}, \gamma_{\mathrm{op}}$, and the number $\nconsider$ of candidate singular vectors to consider.
As a first step, the SVD of the data matrix $\snapshots$ is computed, which is then re-used for fast evaluations of the objective function~\eqref{eq:newlossfast}.
The empty basis $\Vlin$ is initialized alongside index sets $\indexin{0}$ and $\indexout{0}$, that are used to keep track which of the first $\nconsider$ singular vectors are included in the basis.
Then at iteration $i$, a minimizer of \eqref{eq:newlossfast} is determined over the not yet included singular vectors with indices in $\indexout{i}$.
A determined minimizer is then appended to $\scaledVlin_i$ and the index sets are updated, accordingly.
This step is repeated for $\nred$ iterations until we obtain $\scaledVlin = \scaledVlin_r \in \R^{\nfull \times \nred}$.
After $\nred$ iterations, we compute the coefficient matrix $\Vnonlin$ via \eqref{eq:qm_lstsq_problem} and return $\scaledVlin$ and $\Vnonlin$.

\begin{algorithm}[t]
  \caption{Operator inference-aware quadratic manifolds}
  \label{alg:ourquadmani}
  \begin{algorithmic}[1]
    \Procedure{OpInfAwareGreedyQM}{$\snapshots, \nred, \nconsider, \gamma, \gamma_{\bfA}, \gamma_{\bfH}, \gamma_{\mathrm{op}}$}
    \State{Compute the SVD of the snapshot matrix $\leftsings \singvals \rightsings^\top = \snapshots$} 
    \State{Set ${\scaledVlin}_{,0} = [], \indexin{0}=\{\}, \indexout{0}=\{1, \dots, \nconsider\}$}
    \For{$i = 1,\dots, r$}
    \State{Determine $\scaledleftsing{j_i}$ via~\eqref{eq:newlossfast} over all $\{\scaledleftsing{j}\}_{j \in \indexout{i-1}}$ and $\Vnonlin^{\prime} \in \mathbb{R}^{(\nsnapshots -i) \times \nredmod}$}
    \State{Set $\indexin{i}=\{j_1, \dots, j_i\}$ and $\indexout{i}=\{1, \dots, \nsnapshots\}\setminus \indexin{i}$}
    \State{Set $\Vlin_{i}^{\sigma^{-1}} = [\scaledleftsing{j_1}, \dots, \scaledleftsing{j_i}]$}
    \EndFor
    \State{Set $\scaledVlin=[\scaledleftsing{j_1}, \dots, \scaledleftsing{j_r}]$}
    \State{Compute $\Vnonlin$ via the regularized least-squares problem \eqref{eq:qm_lstsq_problem}.}
    \State{Return $\scaledVlin$ and $\Vnonlin$.}
    \EndProcedure
  \end{algorithmic}
\end{algorithm}

\section{Numerical experiments}\label{sec:NumExp}
We now demonstrate the operator inference-aware greedy approach on a transport and turbulent flow example.

\begin{figure}
  \begin{center}
    \resizebox{0.99\textwidth}{!}{\small{\input{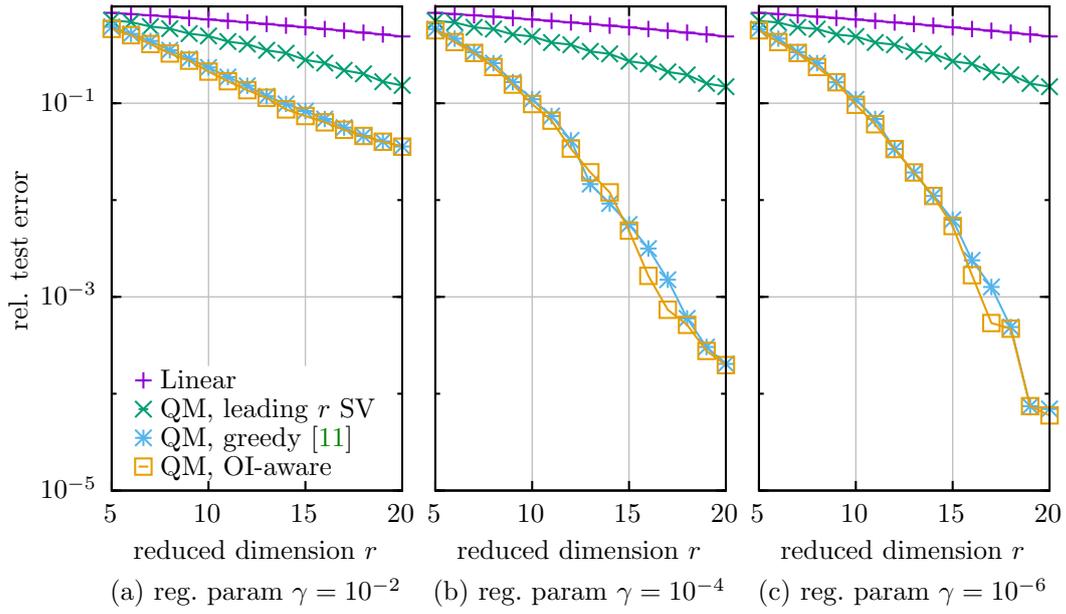}}}
  \end{center}
  \caption{Transport problem: The greedy methods for constructing quadratic manifolds achieve the lowest reconstruction error among all approaches in this example. Furthermore, the error corresponding to the manifold trained with the proposed operator inference-aware (OI-aware) greedy approach is comparable to the error achieved with the regular greedy approach from \cite{SchwerdtnerP2024Greedy}, even though our operator inference-aware approach balances the reconstruction and model prediction errors.}
  \label{fig:transport_projection_err}
\end{figure}

\subsection{Linear transport}
We first consider the linear transport equation. 
\subsubsection{Experimental setup}\label{sec:NumExp:Transport:Setup}
Let $u: [0,T] \times [-\pi,\pi)^2 \to \R, (t,\boldsymbol{x}) \mapsto u(t,\boldsymbol{x})$ be the solution of
\begin{align}\label{eq:NumExp:LinTransport}
  \partial_t u = x_1 \partial_{x_2} u-x_2\partial_{x_1} u
\end{align}
with periodic boundary conditions, where  $\boldsymbol{x}=[x_1, x_2]^\top$.
The initial condition is 
\[
u_0(\boldsymbol{x})=\exp(-80 ((x_1-1)^2+x_2^2)).
\]
To generate snapshot data, we discretize \eqref{eq:NumExp:LinTransport} with 128 spatial Fourier modes in each dimension and apply a fourth-order explicit Runge Kutta scheme with time-step size $\Delta t = 2\pi \times 10^{-3}$.
We simulate the system for $\nsnapshots=2000$ time steps and split the generated data into a training data set 
\[
\snapshots^{\mathrm{(train)}}=[\fullstatei{1},\dots,\fullstatei{1000}]
\]
and a test data set 
\[
\snapshots^{\mathrm{(test)}}=[\fullstatei{1001},\dots,\fullstatei{2000}]\,.
\]
We then train quadratic manifolds with three different techniques. First, we follow \cite{GeelenWW2023Operator,BarnettF2022Quadratic} and use the first $\nred$ left-singular vectors of the training data snapshot matrix. Second, we follow \cite{SchwerdtnerP2024Greedy} and use the greedy method. Third, we use the operator inference-aware greedy method that we introduced in Section~\ref{sec:OpInfGreedy}.
Additionally, we compare to linear approximations using the the reduced space spanned by the first $\nred$ left-singular vectors. 
We consider reduced dimensions $\nred \in \{5,10,15,20\}$. 
The regularization parameter $\gamma$ is determined via grid search over the set $\{10^{-6}, 10^{-4},10^{-2}\}$.

Once we have the encoder-decoder pairs corresponding to either a quadratic manifold or a reduced subspace, we fit an operator inference model to the embedded data as described in Section~\ref{sec:Prelim:OpInf}. The regularization parameters are again determined via grid search over $\gamma_{\bfA}, \gamma_{\bfH} \in \{10^{-4}, 10^{-3},10^{-2},10^{-1}\}$. The value $\gamma_{\mathrm{op}}$ is set to $1.0$ in all experiments.

\begin{figure}
  \begin{center}
    \resizebox{0.99\textwidth}{!}{\small{\input{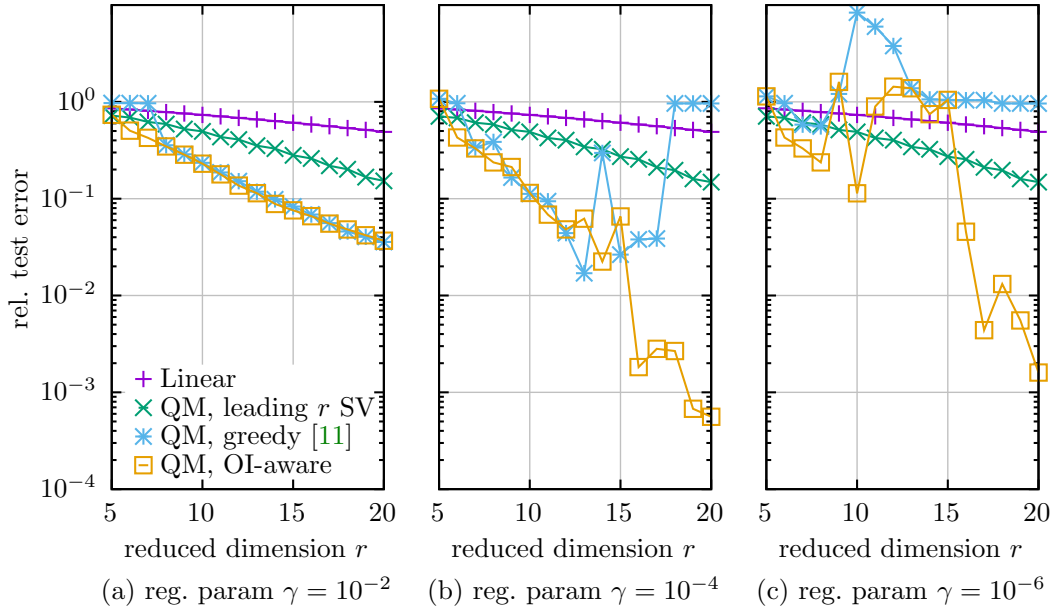}}}
  \end{center}
  \caption{Transport problem: The proposed operator inference-aware (OI-aware) greedy approach leads to operator inference models that achieve  errors that are more than one order of magnitude lower than the error achieved with other reduction methods in this example. }
  \label{fig:transport_opinf_err}
\end{figure}

\begin{figure}
  \begin{center}
    \resizebox{0.99\textwidth}{!}{\small{\input{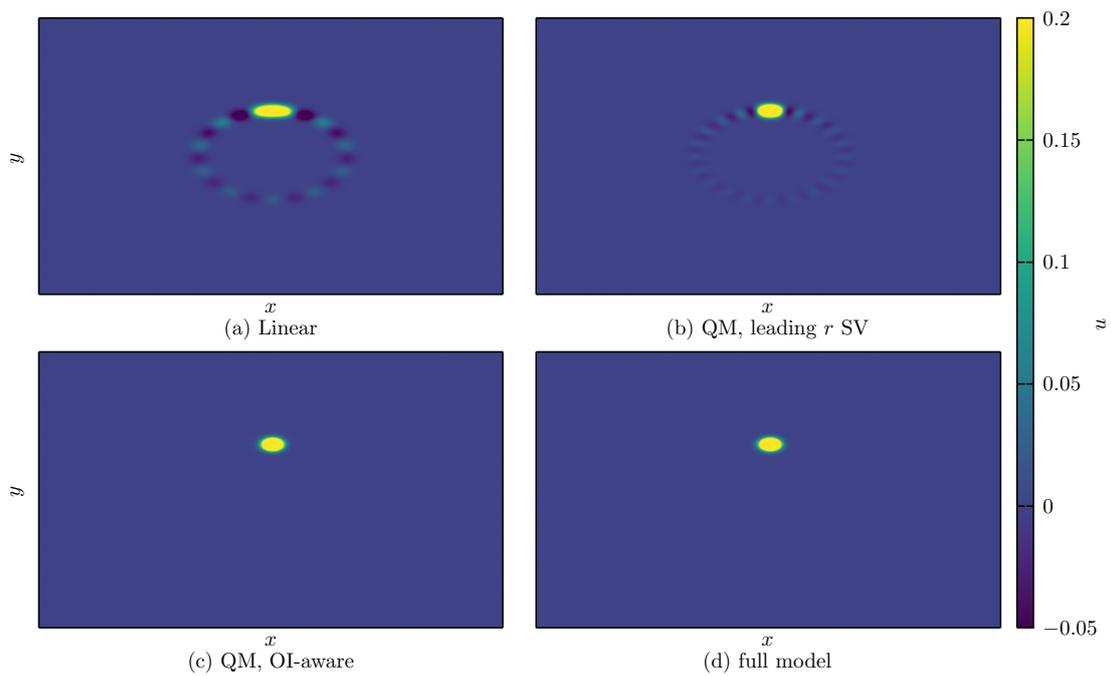}}}
  \end{center}
  \caption{Transport problem: The operator inference models based on linear approximations and on quadratic manifolds trained with the first $\nred$ leading left-singular vectors lead to visible artifacts in the predictions. In contrast, the operator inference model obtained with the proposed operator inference-aware (OI-aware) greedy approach accurately predicts the full-model solution.}
  \label{fig:transport_opinf_reconstruction}
\end{figure}

\subsubsection{Results}
We first consider the relative reconstruction error
\begin{equation}\label{eq:RelReconError}
e_{\mathrm{recon}}(\snapshots^{\mathrm{(test)}})=\frac{\frob{\decoder(\encoder(\snapshots^{\mathrm{(test)}}))-\snapshots^{\mathrm{(test)})}}}{\frob{\snapshots^{\mathrm{(test)}}}},
\end{equation}
where $\decoder$ and $\encoder$ are the decoders and encoders of the different methods we compare. The relative reconstruction error is plotted in Figure~\ref{fig:transport_projection_err}.
First, notice that the results show that decreasing the regularization parameter $\gamma$ from $10^{-2}$ to $10^{-6}$, leads to a decrease in the reconstruction error \eqref{eq:RelReconError}. 
This result is well-aligned with experiments in the literature such as the ones presented in \cite{SchwerdtnerP2024Greedy}. Second, using the leading $\nred$ left-singular vectors only to train the quadratic manifolds leads to a slightly higher accuracy compared to the linear approximation in the reduced subspace, while the greedy approaches achieve orders of magnitude improvements. 
Notice that there is no regularization in case of the reduced subspace and thus the error curve corresponding to the linear approximation is the same across all three plots. 
It is important to note that our new operator inference-aware greedy approach achieves comparable accuracy as the regular greedy approach from  \cite{SchwerdtnerP2024Greedy} in terms of the reconstruction error, even though it additionally uses the model prediction error during the greedy iterations.

Let us now consider the model prediction error of the corresponding operator-inference models,
\begin{equation}\label{eq:MPE}
e_o(\hat{\snapshots}^{\mathrm{(test)}}) = \frac{\frob{\decoder\left( \hat{\snapshots}^{\mathrm{(test)}}\right)-\snapshots^{\mathrm{(test)}}}}{\frob{\snapshots^{\mathrm{(test)}}}},
\end{equation}
where $\hat{\snapshots}^{\mathrm{(test)}}$ is the trajectory predicted for $\encoder(\fullstatei{0,\mathrm{test}})$, where again $\encoder$ and $\decoder$ are the encoders and decoders of the corresponding methods.
Figure~\ref{fig:transport_opinf_err} shows the model prediction error. To keep the visual presentation concise, we again group the plots by the regularization parameter $\gamma$. For the operator-inference models, we show results for the regularization parameters $\gamma_{\bfA}$ and $\gamma_{\bfH}$ with the lowest error. 
The results shown in Figure~\ref{fig:transport_opinf_err} reveal that choosing a higher regularization $\gamma$ can lead to a higher model prediction error. 
However, stronger regularization also means that the  reconstruction error is higher (see Figure~\ref{fig:transport_projection_err}), which is also reflected in the model prediction error. 
For example, for reduced dimension $\nred = 20$, a strong regularization of $\gamma = 10^{-2}$ leads to a model prediction error of about $3 \times 10^{-2}$ for the greedy methods, whereas for the lower regularization $\gamma = 10^{-4}$, the proposed operator inference-aware greedy approach achieves errors as low as $10^{-3}$. 
Furthermore, the results demonstrate that simply applying the regular greedy approach from \cite{SchwerdtnerP2024Greedy} together with operator inference can lead to higher errors, which aligns with the discussion in Section~\ref{sec:Prelim:Problem}.
For example, in Figure~\ref{fig:transport_opinf_err}(b), this is revealed for reduced dimension $\nred=20$, where the error of the operator inference model trained with the quadratic manifold from the regular greedy approach leads to relative errors of about one. 
In contrast, when using our new selection criterion, the model prediction error is as low as $6 \times 10^{-4}$. 
In Figure~\ref{fig:transport_opinf_reconstruction}, we show the prediction of the operator inference model at dimension $\nred = 20$ that achieves the lowest model prediction error for each approach of training the quadratic manifold. 
Artifacts are clearly visible when using a reduced model based on the linear approximation given by the reduced subspace only. Similarly, one can see artifacts when training the quadratic manifold with the leading $\nred$ left-singular vectors only. In contrast, the proposed operator inference-aware greedy approach achieves the lowest error and so avoids these artifacts in the predictions.  

\begin{figure}
  \begin{center}
    \resizebox{0.99\textwidth}{!}{\small{\input{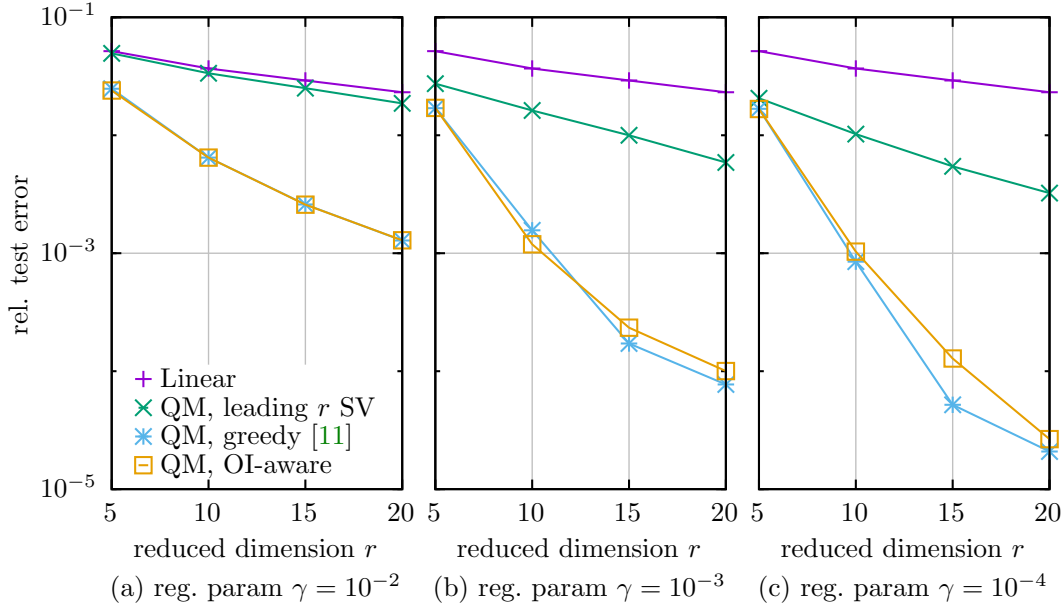}}}
  \end{center}
  \caption{Turbulent flow: The proposed operator inference-aware (OI-aware) greedy approach achieves a comparable reconstruction error as the regular greedy approach even though it balances reconstruction error and model prediction error. Notice that both greedy methods achieve orders of magnitude lower reconstruction errors than the other reduction methods. }
  \label{fig:turbine2d_projection_err}
\end{figure}

\subsection{Flow around a wind turbine}
We now consider data from a large eddy simulation (LES) of an International Energy Agency 15 MW turbine located in a neutral atmospheric boundary layer.

\subsubsection{Experimental setup}
The snapshots for this experiment are generated with AMR-Wind~\cite{Kuhn2025,Sharma2023,sprague2020exawind}, a massively parallel, block-structured adaptive-mesh, incompressible flow solver for wind turbine and wind farm simulations. AMR-Wind solves the incompressible Navier-Stokes equations, as well as temperature, subgrid-scale kinetic energy, and other scalar equations necessary for LES of wind farms. Included physics in these simulations are atmospheric boundary layer forcing, Boussinesq buoyancy, Coriolis forcing, and body forcing to maintain the precursor-derived inflow condition in the presence of the turbine. The turbine model is represented using a Joukowsky disk model~\cite{SORENSEN20202259,Sorenson2023}. The full-model simulation is performed in a $2560m \times 1440m \times 960m$ domain $[x, y, z]^{\top} \in \Omega$, discretized into a grid of $256 \times 144 \times 96$ cells at a time resolution of $0.1s$.
We obtain the dataset after the simulation has reached a fully-developed state by first simulating $600s$ of physical time before data collection to allow the initial transients to dissipate. For this study, we capture hub-height $(z=150m)$ planes (dimension $256 \times 144$) from $k=1000$ consecutive time snapshots, split into training data snapshots
\[
\boldsymbol S^{(\text{train})} = [\fullstatei{0}, \fullstatei{2}, \dots, \fullstatei{998}] \in \R^{36864\times 500}
\]
and test data snapshots
\[
\boldsymbol S^{(\text{test})} = [\fullstatei{1}, \fullstatei{3}, \dots, \fullstatei{999}]\in \R^{36864\times 500}\,.
\]
Training the quadratic manifolds and operator-inference models is analogous to the previous example as described in Section~\ref{sec:NumExp:Transport:Setup} with
the grid search for the regularization parameter $\gamma$ being over the set $\{10^{-4}, 10^{-3},10^{-2}\}$ and $\gamma_{\mathrm{op}}$ set to $10^2$ in all experiments.

\begin{figure}
  \begin{center}
    \resizebox{0.99\textwidth}{!}{\small{\input{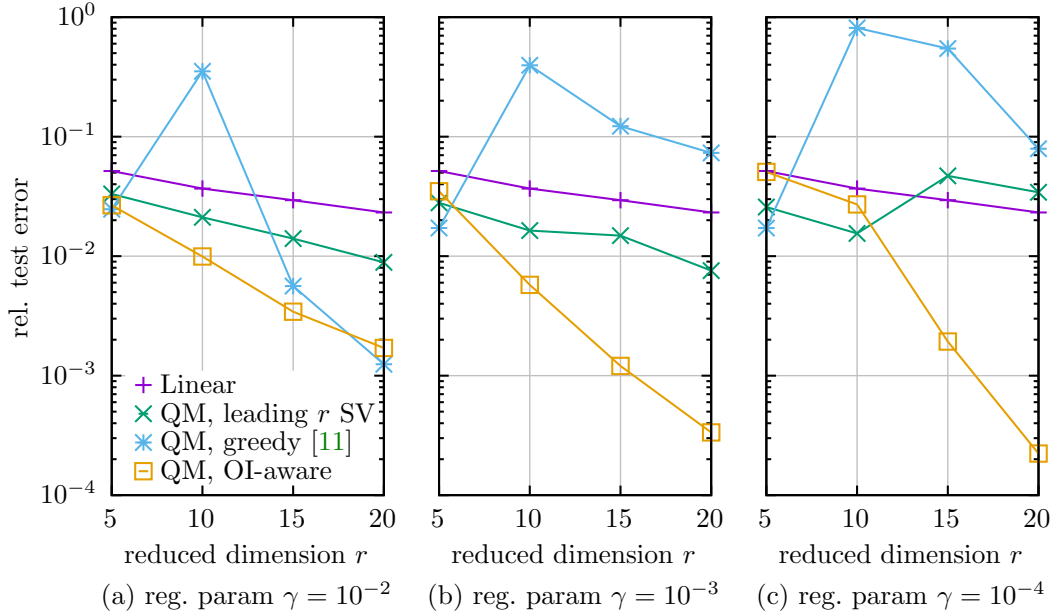}}}
  \end{center}
  \caption{Turbulent flow: Operator inference models based on the proposed operator-inference aware (OI-aware) greedy approach achieve orders of magnitude lower errors than models learned with other reduction methods in this example. }\label{fig:turbine2d_opinf_err}
\end{figure}

\begin{figure}
  \begin{center}
    \resizebox{0.99\textwidth}{!}{\small{\input{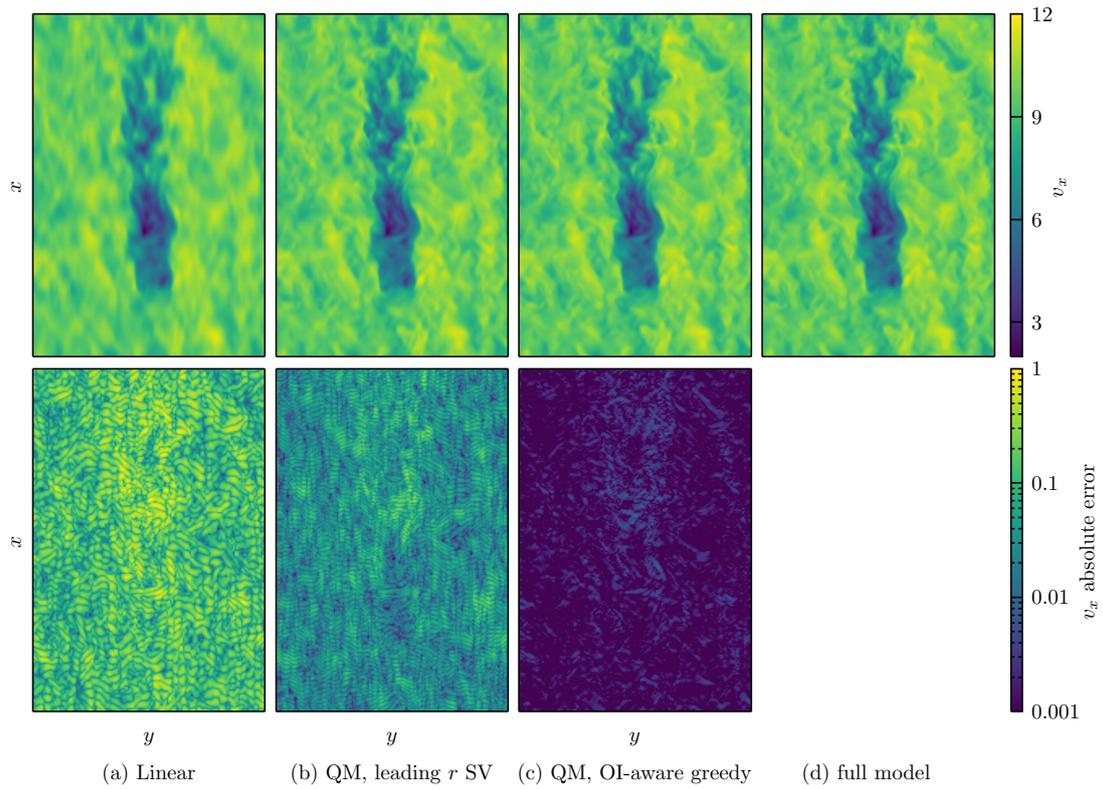}}}
  \end{center}
  \caption{Turbulent flow: The errors shown in the bottom row demonstrate that the proposed operator inference-aware greedy approach achieves the lowest error in this example. }\label{fig:turbine2d_opinf_reconstruction}
\end{figure}

\subsubsection{Results}

In Figure~\ref{fig:turbine2d_projection_err}, we study the relative reconstruction error \eqref{eq:RelReconError}. Recall that the reconstruction error only captures the error made by the encoder-decoder pair and does not involve the operator-inference model. 
The greedy methods lead to a error reduction of up to two orders of magnitude compared to the quadratic manifold constructed from the first $\nred$ leading left-singular vectors .
The operator inference-aware greedy approach in this work avoids a loss of accuracy compared to the regular greedy from \cite{SchwerdtnerP2024Greedy} even though the additional criterion based on the model prediction error is included.

Let us now consider the model prediction error \eqref{eq:MPE}, which is reported in Figure~\ref{fig:turbine2d_opinf_err}. 
The proposed operator inference-aware greedy method achieves the lowest error on the test data for $\gamma = 10^{-3}$ and $\gamma = 10^{-4}$. 
In contrast, the regular greedy approach fails to provide meaningful predictions in this setting for $\gamma = 10^{-3}$ and $\gamma = 10^{-4}$ with relative errors higher than 10\%. 
We conclude that these results show that including the model prediction error into the greedy selection leads to quadratic manifolds that are better suited for operator inference in the sense that orders of magnitude lower prediction errors are achieved. 

\begin{figure}
    \centering
    \resizebox{0.99\textwidth}{!}{\small{\input{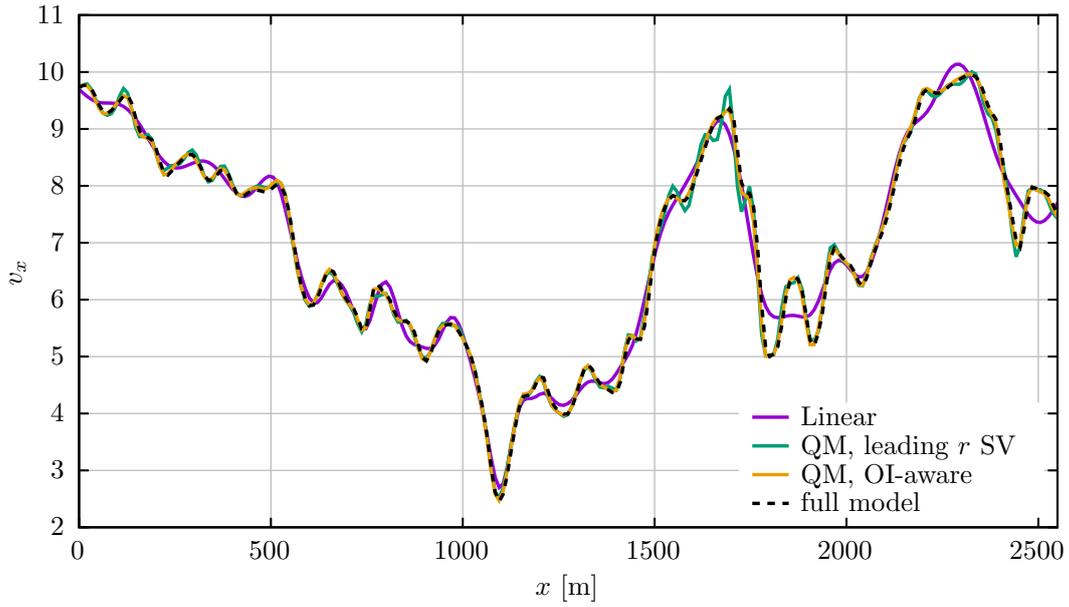}}}
    \caption{Turbulent flow: The prediction of the wind velocity $v_{x}$ after the turbine at $x = 720m$ shows that operator inference models based on linear approximations lead to a smoothing in the prediction that misses the oscillatory behavior of the quantity of interest. The model trained on the quadratic manifold based on the first $\nred$ leading left-singular vectors alone leads to overshooting in the prediction of the quantity of interest, e.g., near $x = 1600m$. In contrast, the models based on the proposed operator inference-aware (OI-aware) greedy lead to accurate predictions of the quantity of interest over the whole $x$ range of interest.}
    \label{fig:qoi}
\end{figure}

Reconstructions of the final-time  state using  the setting with the lowest model prediction error for each method are shown in Figure~\ref{fig:turbine2d_opinf_reconstruction} for dimension $\nred=20$. First notice that the operator inference model based on the linear approximation in the reduced space spanned by the first $\nred$ left-singular vectors leads to high errors. Second, the quadratic manifold trained with the first $\nred$ leading left-singular vectors leads to lower errors than linear approximations. However, the lowest error is obtained with the operator inference model corresponding to the quadratic manifold trained with the proposed operator inference-aware greedy method. 
In Figure~\ref{fig:qoi}, we show the flow downstream of the wind turbine at $x=720m$, which is a quantity of interest that is essential for computing power generation. The proposed greedy method leads to an operator inference model that predicts the quantity of interest with highest accuracy among the methods that are compared here.  

\section{Conclusions}\label{sec:Conc}
The proposed operator inference-aware greedy approach demonstrates orders of magnitude higher accuracy can be achieved when the training of encoder-decoder pairs is not based only on the reconstruction error on snapshot data but explicitly takes into account the model prediction error. In the specific case that we consider here, namely quadratic manifolds trained with greedy methods, the results show that the proposed operator inference-aware greedy avoids embeddings that are oscillatory or otherwise non-smooth. Avoiding non-smooth embeddings leads to operator inference models that achieve up to two orders of magnitude higher accuracy than models based on quadratic manifolds that are trained with the regular greedy approach that ignores model prediction errors in the selection.

\section*{Declarations}
\textbf{Funding statement}
The authors were supported by the US Department of Energy, Office of Scientific Computing Research,
 DOE Award DE-SC0024721. 
(Program Manager Dr. Margaret Lentz).
This work was authored in part by the National Renewable Energy Laboratory for the U.S. Department of Energy (DOE) under Contract No. DE-AC36-08GO28308. A portion of the research was performed using computational resources sponsored by the DOE Office of Energy Efficiency and Renewable Energy and located at the National Renewable Energy Laboratory. The views expressed in the article do not necessarily represent the views of the DOE or the U.S. Government. The U.S. Government retains and the publisher, by accepting the article for publication, acknowledges that the U.S. Government retains a nonexclusive, paid-up, irrevocable, worldwide license to publish or reproduce the published form of this work, or allow others to do so, for U.S. Government purposes. 

\textbf{Competing interests} The authors declare no competing interests.

\textbf{Author Contribution}
Paul Schwerdtner:
Conceptualization; Methodology; Formal analysis and investigation; Writing - original draft preparation; Writing - review and editing;
Prakash Mohan:
Conceptualization; Formal analysis and investigation; Writing - original draft preparation; Writing - review and editing;
Julie Bessac:
Conceptualization; Writing - original draft preparation; Writing - review and editing; Funding acquisition; Resources; Supervision;
Marc T. Henry de Frahan:
Conceptualization; Writing - original draft preparation; Writing - review and editing;
Benjamin Peherstorfer:
Conceptualization; Writing - original draft preparation; Writing - review and editing; Funding acquisition; Resources; Supervision;

\bibliographystyle{my-sn-mathphys-num}
\bibliography{donotchange,references}

\end{document}